\documentclass[a4paper]{article}
\usepackage{natbib}
\pdfoutput = 1

\usepackage{graphics,graphicx}
\usepackage{amsmath}
\usepackage{amsfonts}
\usepackage{float,times,amssymb,bm}
\usepackage{colortbl}

\bmdefine\bX{X}
\bmdefine\bTheta{\Theta}

\newcommand{\xihat}{\ensuremath{\hat{\xi}}}
\newcommand{\psihat}{\ensuremath{\hat{\psi}}}
\newcommand{\asinh}{\ensuremath{\text{asinh}}}

\begin{document}

\title{A Noninformative Bayes-like Approach to Probability-Preserving Prediction of Extremes.}
\author{Allan McRobie \\
Cambridge University Engineering Department\\
Trumpington St, Cambridge, CB2 1PZ, UK \\
fam20@cam.ac.uk}

\maketitle

\begin{abstract}
The extrapolation of extremes to values beyond the span of stationary univariate historical data is considered from Bayesian and Frequentist perspectives.
The intention is to make predictions which in some sense ``preserve probability''.
A Frequentist approach based on a simple curve-fit estimate of the tail parameter $\xi$ of a Generalised Pareto Distribution was described in McRobie (2014) (arXiv:1408.1532). In this paper, the corresponding Bayes-like approach is described, using a plausible noninformative prior for the tail parameter. The two approaches, though philosophically different, show a reasonable degree of correspondence.
\end{abstract}

\section{Introduction}
\cite{McRobiecurvefit} presented a method for extrapolating to extreme values outside the span of historical data.
The algorithm was based on sampling theory applied to the tail of the Generalized Pareto Distribution (GPD), and led to a location- and scale-invariant predictor which preserved probability to a very good approximation. The paper here considers that same approach, but from a Bayes-like perspective, using a non-informative prior on the GPD tail parameter.

Given that this paper extends the results of a series of papers by the author on this topic, we do not provide here a full literature review of the wider prediction problem, and refer only to the papers in the preceding series.

For a process generating samples with an underlying stationary probability distribution $F(x \ | \ \bTheta^*)$
there is a true value $x_{T,true}$ of the $T$-level event, this being such that $G(x_{T,true} \ | \ \bTheta^*) = 1- F(x_{T,true} \ | \ \bTheta^*) = 1/T$. That is, the probability that the next singleton drawn exceeds $x_{T,true}$ is $1/T$.
Traditional approaches to extreme value analysis aim to provide an estimate $x_{T,est}$ which is in some sense close to $x_{T,true}$,
usually via first making an estimate $\hat{\bTheta}$ of the true but unknown parameter set  $\bTheta^*$. The probability that the next sample drawn will exceed $x_{T,est}$ may however be very different to $1/T$. This is explored in \cite{McRobieEEE}.

The Bayesian framework provides an alternative perspective on this problem. Rather than {\it estimating} the $T$-level event, the Bayesian creates a {\it predictive} distribution, this being a distribution of beliefs (given the data and the prior knowledge) about the next
singleton drawn. The Bayesian $T$-level prediction $x^B_{T,pred}$ is then that value of $x_{next}$ for which a proportion $1/T$ of the Bayesian's beliefs lie above $x^B_{T,pred}$. It has a $T$-level exceedance probability in the sense that a Bayesian would consider the possibility that the next singleton drawn will exceed their prediction to be as surprising as winning a single bet on a fair wheel-of-fortune with $T$ equal divisions.

Recognising that an estimate of $x_{T,true}$ may be exceeded more (or less) often than once in every $T$ events, a Frequentist may instead choose to adopt a more unconventional strategy, and - somewhat emulating the Bayesian - make a prediction $x^F_{T,pred}$ with the property that, under repeated sampling and prediction at any fixed parameter set $\bTheta^*$, the next samples drawn exceed their corresponding predictions a proportion $1/T$ of the times. A predictor with such a property is here called a ``Probability Preserving Predictor'' (and was called an ``Exact Exceedance Estimator'' in \cite{McRobieEEE}).

\begin{figure} [h!] \centering
   \includegraphics[width=110mm, keepaspectratio]{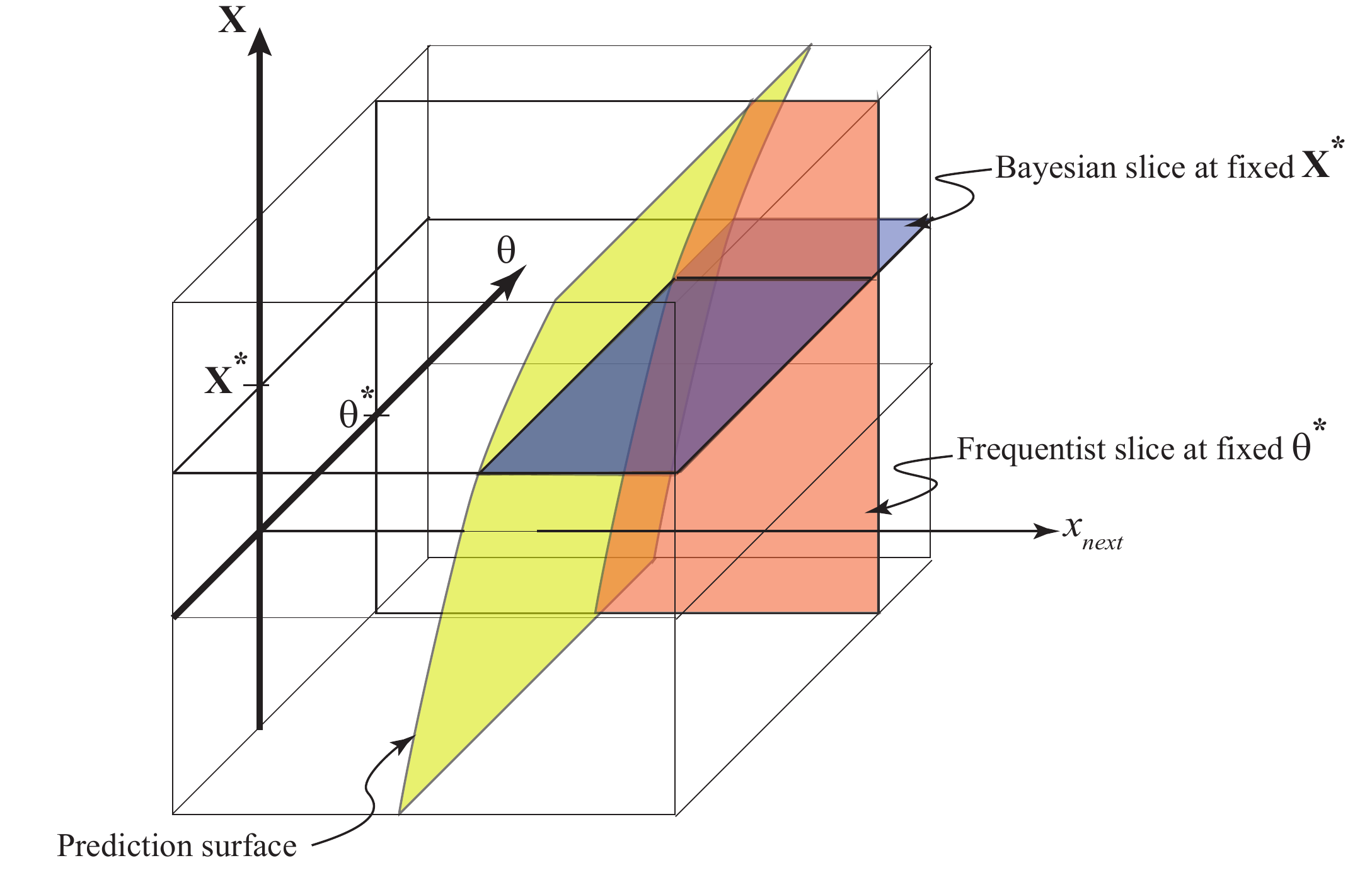}
 \caption{ A schematic of the full problem space $\lbrace \bX, \bTheta, x_{next} \rbrace$, with the subspaces of data $\bX$, parameters $\bTheta$ and next sample $x_{next}$.  Frequentists focus attention on the vertical slices $\lbrace \bX, x_{next} \ | \ \bTheta^* \rbrace$ whilst Bayesians consider
 horizontal slices $\lbrace \bTheta, x_{next} \ | \ \bX^* \rbrace$.}
\label{predictionsurf}
\end{figure}

Figure~\ref{predictionsurf} shows the prediction problem seen from the two different philosophical perspectives.
The space is spanned by the data $\lbrace \bX \rbrace$ (of dimension $N$), the parameters $\lbrace \bTheta \rbrace $ (of dimension $m$) and the next singleton sample drawn $\lbrace x_{next} \rbrace$ (of dimension 1). At desired recurrence level $T$, the prediction surface $x_{T,pred}(\bX)$ is an $N$-dimensional surface, a function over the data space $\lbrace \bX \rbrace$. The data space $\lbrace \bX \rbrace$ and the parameter space $\lbrace \bTheta \rbrace $ axes are each shown by thicker arrows to denote that each is typically multi-dimensional. Frequentists and Bayesians may not agree to construct the same prediction surface, but only one surface is drawn in the Figure to illustrate the differing ways that the two analysts approach the problem.

For the Frequentist, attention is focused on an $(N+1)$-dimensional vertical slice at the actual parameter vector $\bTheta^* \in \lbrace \bTheta \rbrace$. At any given $\bTheta^* \in \lbrace \bTheta \rbrace$, there is a true value $x_{T,true}$ such that $P(x_{next}>x_{T,true} \ | \ \bTheta^*) = 1/T$, irrespective of the data. As described above, an orthodox Frequentist will, via estimates of the parameter, endeavour to construct an estimate $x_{T,est}$ of $x_{T,true}$, perhaps even constructing confidence intervals around $x_{T,est}$. The Frequentist approach adopted in this paper is completely orthodox, excepting only that, rather than constructing the estimate $x_{T,est}$, the Frequentist is willing to endeavour to construct a prediction surface $x_{T,pred}^F(\bX)$ such that the probability that $x_{next} > x_{T,pred}^F$ is $1/T$ for any fixed $\Theta^*$.

For the Bayesian, attention focuses instead on a $(m+1)$-dimensional horizontal slice at the given data $\bX^* \in  \lbrace \bX \rbrace$. Unlike the Frequentist, the Bayesian entertains the notion of ``a probability of the parameters''. This begins with a distribution of prior beliefs $\Pi_0(\bTheta)$ over the parameter space $\lbrace \bTheta \rbrace $, and is then conditioned on the given data  $\bX^*$ to obtain the posterior belief distribution $\Pi_1(\bTheta \ | \ \bX^*)$ over parameters. The Bayesian then constructs the predictive distribution
\begin{equation}
p(x_{next} \ | \ \bX^*) = \int_{\forall \bTheta} p(x_{next} \ | \ \bTheta) \  \Pi_1(\bTheta \ | \ \bX^*) \ d\bTheta
\end{equation}
which represents the analyst's beliefs about the next singleton (given the data and the prior beliefs about parameters). The Bayesian's $T$-level prediction $x_{T,pred}^B$ is then determined by the tail integral of the predictive distribution, being that value for which
\begin{equation}
\frac{1}{T} =  \int_{x_{T,pred}^B}^\infty  p(x_{next} \ | \ \bX^*) \ dx_{next}
\end{equation}

The Bayesian, thus far, is completely orthodox, and the prior $\Pi_0(\Theta)$ would be a proper distribution that appropriately represents their subjective beliefs regarding the parameters before they have looked at the data. The first step in relaxing this orthodoxy here is to allow the Bayesian to hold improper distributions of prior beliefs. These may include beliefs that are uniform over the whole real line for location parameters and the well-known $1/\sigma$ prior for scale parameters. An earlier paper (\cite{McRobieEEE}) demonstrated that for prediction problems on any location-and-scale distribution $F((x-\mu)/\sigma)$, a Bayesian predictor constructed from a $1/\sigma$ prior is necessarily probability-preserving for samples drawn at ANY parameter $\bTheta^* \in \lbrace \bTheta \rbrace $. In that case, then, there is a rather beautiful correspondence between the Bayesian and Frequentist approaches, and both schools could agree on a prediction surface $x_{T,pred}(\bX)$.

The key point is that, to construct a prediction surface, Frequentists perform integrals over vertical slices (at $\bTheta^*$ fixed) whilst Bayesians perform integrals over horizontal slices (at $\bX^*$ fixed). Although the domains of integration are typically distant and disparate regions of the full $\lbrace \bX, \bTheta, x_{next} \rbrace$ problem space, it was shown in \cite{McRobieEEE} that for location-and-scale problems there is a natural transformation between the two regions such that the integrals turn out to be equivalent.  This is somewhat remarkable, given that the two integrals are not only over different regions of the problem space but are generally over spaces of differing dimensions, the Bayesian integral over parameters being 2-dimensional in the location-scale case, and the Frequentist integral being over the $N$-dimensional data space. (The identification is achieved by partitioning the Frequentist's $N$-dimensional data space into a collection of slices, each slice containing a 2D plane on which probability densities may be mapped across to those in the Bayesian's perspective).

The intention of the investigation here is to see if the pleasing correspondence between the Frequentist and $1/\sigma$-Bayesian perspectives found in the two-parameter $(\mu, \sigma)$ location-and-scale case can be generalised to the three-parameter $(\mu, \sigma, \xi)$ case of the Generalized Pareto Distribution (GPD) for use in extreme value analysis. It is not obvious that there should be any such neat correspondence. However, \cite{McRobiecurvefit} demonstrated the existence of a simple location-and-scale invariant predictor which was based on a Frequentist perspective and gave very close to probability-preserving performance over all tail parameters $\xi$ for the GPD. This leads naturally to the question as to whether there exists a nearby Bayesian predictor which gives similar results. The question is complicated considerably by the addition of the third parameter $\xi$. The location-and-scale case benefits from the properties of affine transformations, and these are lost with the addition of the third parameter. Moreover, many of the integrals become extremely complicated analytically. For example, \cite{McRobiePRED} made only limited progress in constructing an (approximately) probability-preserving predictor for the GPD due to the complexity of the numerous hypergeometric and Lauricella functions that were encountered in the analysis. These functions, and the probability densities that they represent, are replete with numerous singularities, making the analysis rather tricky. The curve-fit approach adopted in  \cite{McRobiecurvefit}, however, was considerably simpler analytically and led to a rather natural procedure for extrapolating to extreme values beyond the span of historical data. Although the construction was based on Frequentist sampling theory, the construction of a predictor $x_{T,pred}^F (\bX)$ at any recurrence level $T$ implies the existence of a form of  ``predictive  distribution'' highly reminiscent of the Bayesian predictive distribution. The objective here, then, is to construct a Bayes-like approach (based on integrals over horizontal slices)  that in some sense matches the Frequentist curve-fit approach (which was based on integrals over vertical slices).

The approach is called ``Bayes-like'' rather than Bayesian, because we shall need to ask the Bayesian to be a little more unconventional than merely accepting the possibility of improper priors. The scenario adopted here is that the Bayesian, rather than being given access to all the data $\bX$, will only be given access to the curve-fit estimate $\hat{\xi}(\bX)$. (Strictly speaking, the Bayesian will also eventually need knowledge of the two data points used in the location-scale normalisation). Moreover, integrals over the parameter space will be restricted to one dimensional integrals involving only the tail-parameter $\xi$, as a result of the removal of the location- and scale-information by restricting consideration to the location- and scale-invariant statistic $\hat{\xi}$.

\begin{figure} [h!] \centering
   \includegraphics[width=110mm, keepaspectratio]{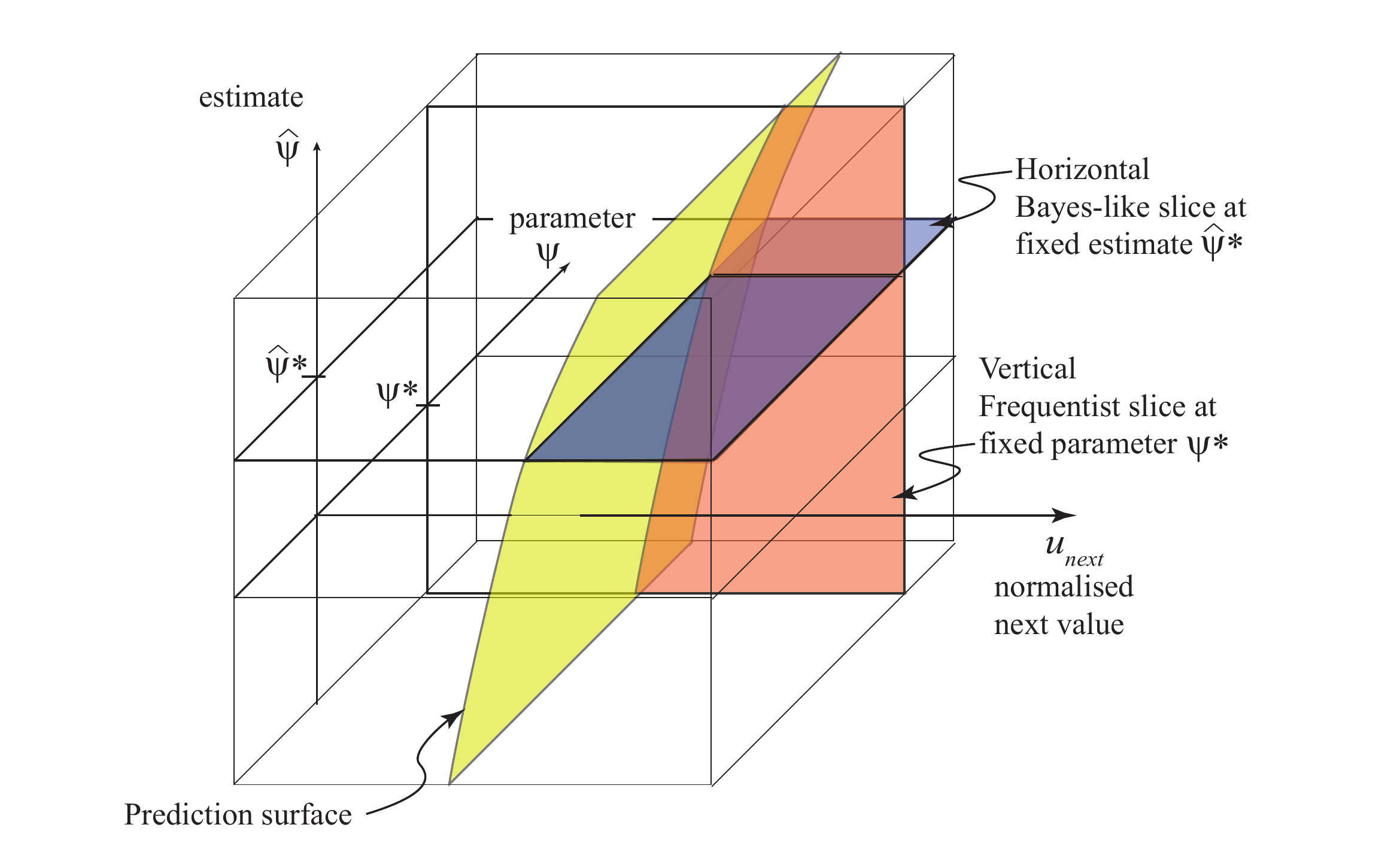}
 \caption{The reduced (3D) problem space that will be adopted in this paper. By normalising the data to remove location and scale effects, the parameter space is one dimensional, with $\psi$ being some function of the tail parameter $\xi$. The data space $\bX$ has been projected down to the single dimension, a statistic $\psihat(\bX)$ (which is an estimator of $\psi$). The next singleton drawn - although now normalised - again defines a one-dimensional subspace.}
\label{reducedpredictionsurf}
\end{figure}

The general prediction scheme illustrated earlier in Figure~\ref{predictionsurf} has thus now evolved to the specific scheme illustrated in Figure~\ref{reducedpredictionsurf}. The parameter space $\lbrace \psi \rbrace$ and the data space $\lbrace \hat{\psi} \rbrace$ are now each one dimensional, and the parameter $\psi$ is merely a re-parameterisation of the GPD tail parameter $\xi$ such that the prior distribution $\Pi_0$ on $\psi$ is improper uniform over the whole real line.  That is, prior beliefs are represented by an appropriate re-parameterisation of the form $\psi = f(\xi)$ such that $\Pi_0(\psi) = 1 $ for $-\infty < \psi < \infty$. The choice of the function $f$  as $\psi = f(\xi) = \asinh(\xi)$ will be described in the following section. The data space, likewise, has been reduced down to the single dimension $\lbrace \hat{\psi} \rbrace$, where $\hat{\psi} = f(\hat{\xi}) = \asinh(\hat{\xi})$, with $\hat{\xi}$ the curve-fit estimate of $\xi$. As a result of the normalisation to a location- and scale-invariant problem, the single-dimensional space representing the next singleton drawn is now no longer simply $\lbrace x_{next} \rbrace$, but the normalised $u_{next} = (x_{next}- x_{N/2})/(x_{N/2} - x_{N})$. This will be described in the next section. The normalised next singleton drawn thus depends on the data via the two data points $x_{N/2}$ and $x_{N}$, which is a complication. Furthermore, we choose to reparameterise the normalised next singleton drawn as $w_{next} = \asinh(u_{next})$. This is merely for computational convenience in order to have some control over possibly large numerical values. Given that we are only concerned with the tail cumulative distribution $G(u_{next})$, this transformation is of no mathematical or philosophical significance, because cumulative distributions are invariant under monotonic transformation.

The final nuance in the Bayes-like construction is that the density on the full space no longer decomposes as neatly as it does in the pure (Figure~\ref{predictionsurf}) case. There, the overall density over $\lbrace \bX, \bTheta, x_{next} \rbrace $ could be decomposed into a product of the posterior density $\Pi_1(\bTheta | \bX) $  on $\bTheta \times \bX$ subspaces and the basic density $p(x_{next} | \bTheta)$ on $\bTheta \times x_{next}$ subspaces. Although in the Bayes-like construction, the posterior simplifies to the likelihood function via the use of the improper uniform prior, the behaviour in the ``next sample drawn'' direction (i.e. $u_{next}$) is now more complicated, since $u_{next}$ has data dependence via the $x_{N/2}$ and $x_{N}$ terms used in the location-scale normalisation. Considering the projection down from the full $(N+m+1)$-dimensional space $\lbrace \bX, \bTheta, x_{next} \rbrace $ to the 3-dimensional $\lbrace \psihat, \psi, u_{next} \rbrace $ space, the resulting density on the smaller space is
\begin{equation}
p(\psihat, \psi, u_{next})  = \int  p(\bX,\bTheta,x_{next}) \frac{\partial (\bX, \bTheta, x_{next} ) }{ \partial (\psihat, \psi, u_{next}) } \ dV
\end{equation}
where the integral is over all unwanted variables, (namely $\mu$ and $\sigma$ in the parameter space, and any $N-1$ data variables that are transverse to surfaces of constant $\psihat$ in the data space). Although in theory possible, the analytical evaluation of the right-hand side is extremely complicated. However, as will be demonstrated in later sections, it can be readily approximated numerically.

The resulting predictive density can only be decomposed as far as
\begin{equation}
p(\psihat, \psi, u_{next})  = p(u_{next} | \psihat, \psi) p(\psihat, \psi)
\end{equation}
where the latter factor, for improper uniform prior on $\psi$, is simply the likelihood function. The point is that, in the ``next sample'' direction, we can no longer use the simple basic (data independent) distribution, but must use the more general density $ p(u_{next} | \psihat, \psi)$ which depends on both the parameter $\psi$ AND the data (via $\psihat$).

We are thus asking the Bayesian to accept the use of improper priors, and to be willing to make predictions in the case where they are only given access to the statistic $\psihat$. We thus call the approach ``Bayes-like'', rather than Bayesian, but note that it still entertains the notion of ``a probability of a parameter'' and it still constructs the predictor by performing integrals over horizontal slices at fixed ``data'' $\psihat$, in contrast to the Frequentist integrals over vertical slices at fixed parameter $\psi$.

\section{A candidate noninformative prior for the tail parameter}
As described in \cite{McRobiecurvefit}, the location- and scale-invariant curve-fit estimator $\xihat$ of the tail parameter $\xi$ of GPD data is constructed by fitting a curve through upper order statistics. The data is ordered such that $x_1$ is the largest upper order statistic and $x_N$ the smallest.  Only the case $N=20$ was considered, although generalisation to any $N$ is straightforward.  For this case, the procedure fits a curve through the upper nine order statistics, normalised with respect to the tenth and twentieth: \begin{equation}
u_i = \frac{x_i - x_{10}}{x_{10} - x_{20}} \ \ \  \text{    with } 1 \leq i \leq 9
\label{scaleddata}
\end{equation}
and the curve fitting uses the logarithm $\log(1+u_i)$,  selecting the estimate $\xihat$ as that which minimises the sum of the squares of the nine residuals
\begin{equation}
\epsilon_i  = \log(1+u_i) - \log(1+u) \ \ \text{     with    }  u = \frac{(G_{10}/G_i)^{\xihat} -1 }{1 - (G_{10}/G_{20})^{\xihat}}
\end{equation}
the tail exceedance probabilities $G_j$ being approximated as $(j-0.5)/N$ in this estimation phase.
This is little other than estimating $\xi \approx \xihat$ via a curve-fit to the empirical distribution of the normalised data.

\begin{figure} [h!] \centering
   \includegraphics[width=72mm, keepaspectratio]{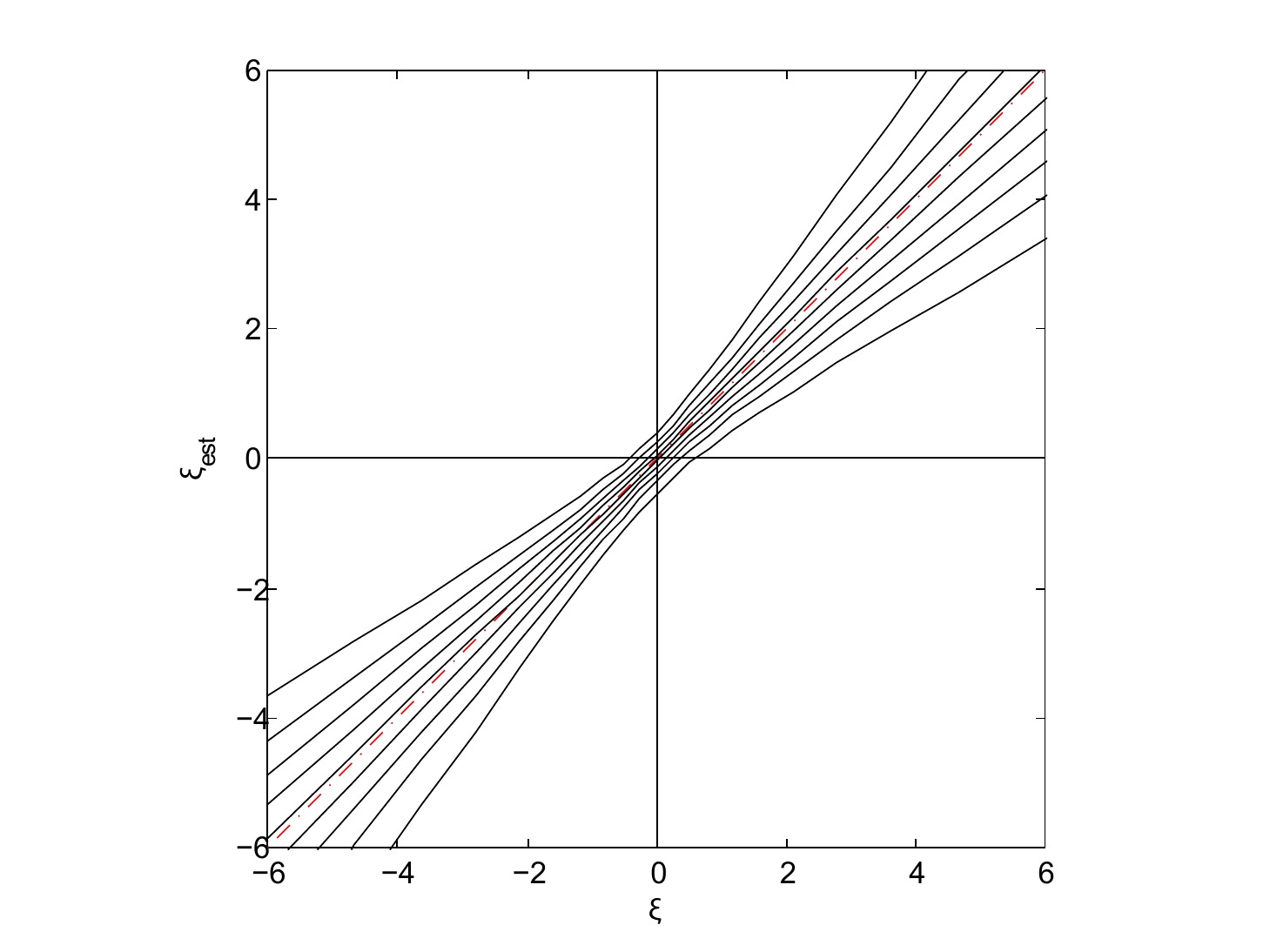}
   {\includegraphics[width=72mm, keepaspectratio]{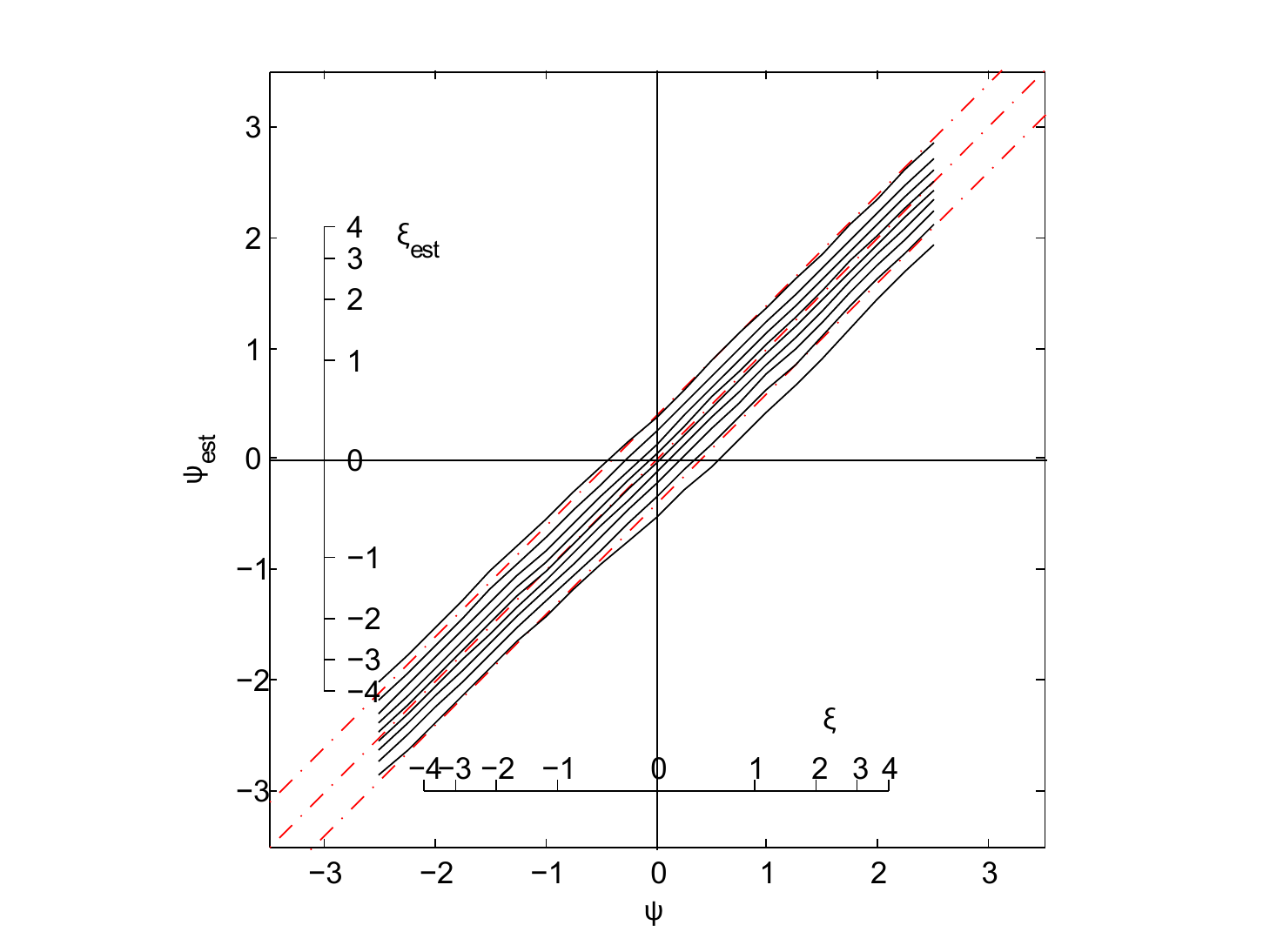}}\\
   {\it a}) \hspace{72mm} {\it b}) \hspace{72mm}
 \caption{Deciles of the likelihood of the estimator, plotted in {\it a}) as $\xihat$ versus $\xi$, and in {\it b}) as $\psihat$ versus $\psi$, where $\psi = \asinh(\xi)$ and $\psihat = \asinh(\xihat)$.}
\label{Psilikelihood}
\end{figure}

\begin{figure} [ht!] \centering
   \includegraphics[width=70mm, keepaspectratio]{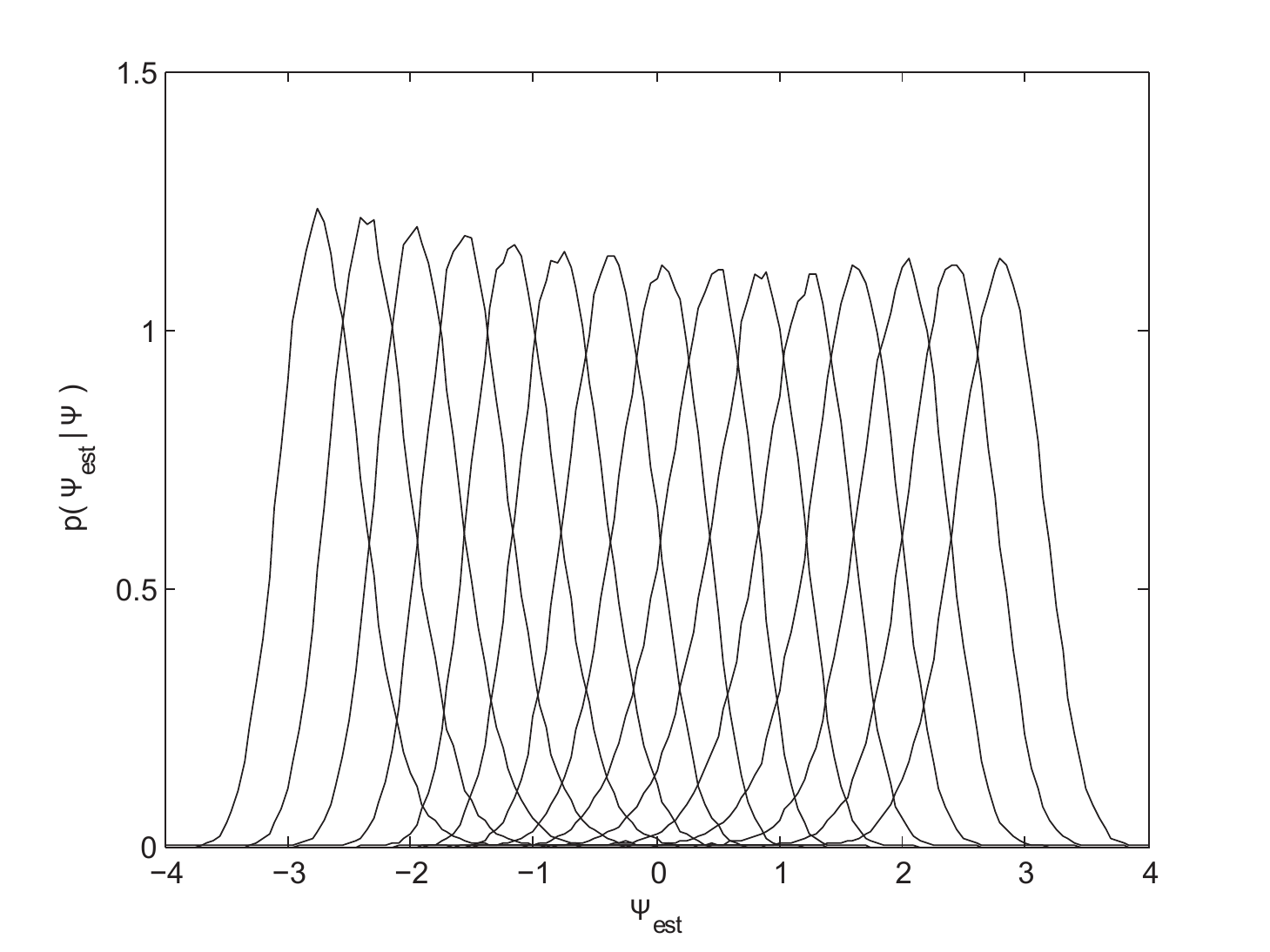}
   \includegraphics[width=70mm, keepaspectratio]{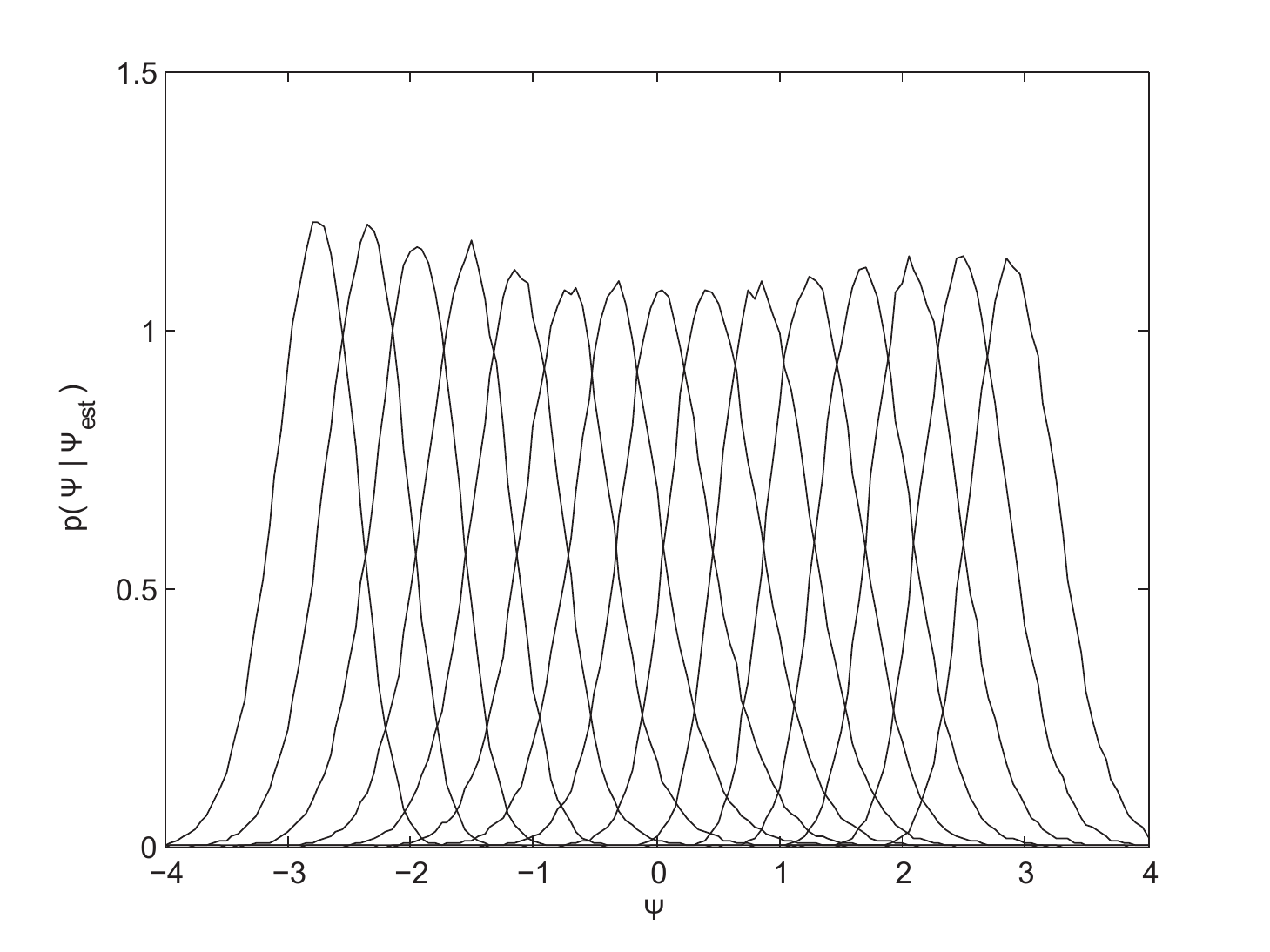}\\
    \includegraphics[width=70mm, keepaspectratio]{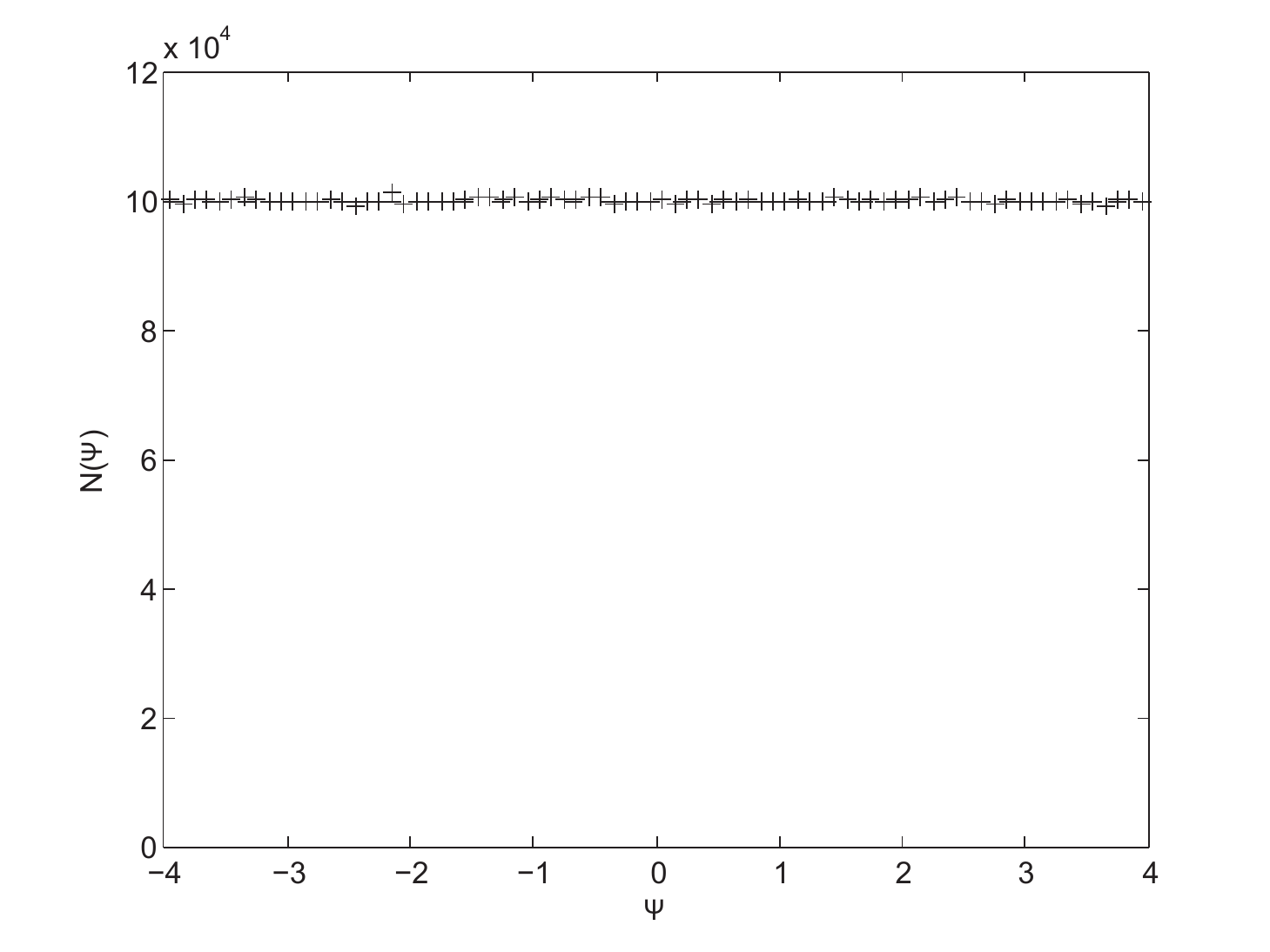}
   \includegraphics[width=70mm, keepaspectratio]{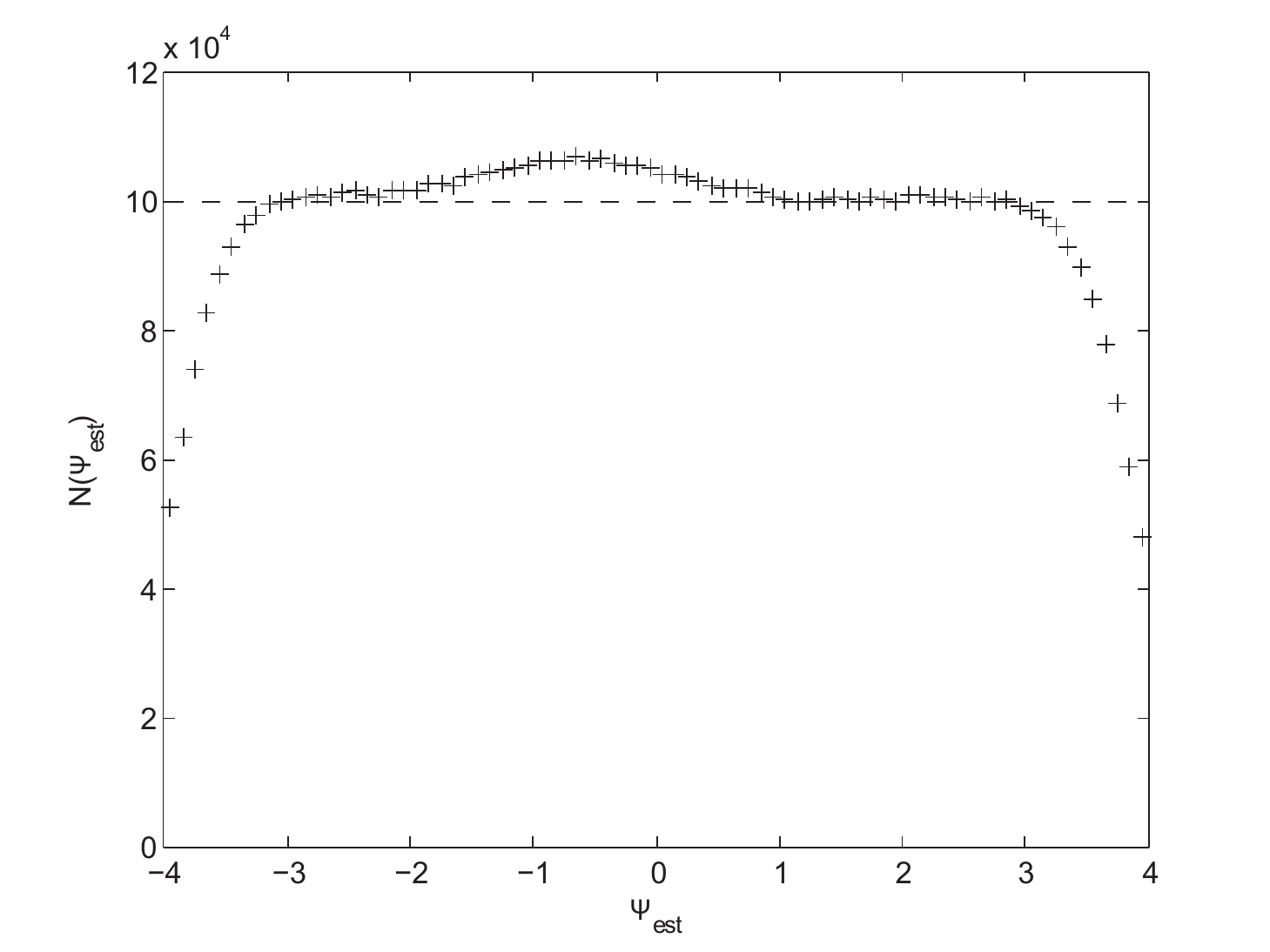}
 \caption{The upper diagrams show (normalised) densities on slices through the likelihood function, at constant $\psi$ (left) and constant $\psihat$ (right). The left-hand figures correspond to the standard Frequentist sampling density of the estimator, and the right-hand figures correspond to the Bayesian posterior, ``the probability of the parameter, given the data''.  (To avoid clutter, only every fourth density computed has been plotted). The lower figures show the number of data points in each slice. Given that the parameter $\psi$ was uniform over $[-4,4]$, the number in each vertical Frequentist slice was approximately constant. The horizontal Bayesian slices, however, show end effects as a result of some estimates lying outside the range $[-4,4]$. To avoid the influence of such end effects, Bayesian posteriors were thus only constructed for estimates $\psihat$ in the range $[-3,3]$.}
\label{psigivenpsihat}
\end{figure}

Figure~\ref{Psilikelihood}{\it a}  illustrates the distribution of the resulting estimates $\xihat$ as the underlying parameter $\xi$ is varied. It shows that the variance of the estimate $\xihat$ is minimal near $\xi = 0$,  and grows as the magnitude of $\xi$ increases away from this central region. These same results are plotted in Figure~\ref{Psilikelihood}b) after the re-parameterisations $\psi = \asinh(\xi)$ and $\psihat = \asinh(\xihat)$. In this re-parameterisation, the variance of the estimate $\psihat$ remains almost constant with parameter $\psi$. Indeed, the resulting distribution (which is essentially the likelihood function) is almost location-invariant along the $\psi = \psihat$ diagonal. This can be also be seen in Figure~\ref{psigivenpsihat} later, which shows slices through the distribution. Given that the $1/\sigma$ prior for location- and scale-invariant problems reduces to the improper uniform prior for location-invariant problems, the almost-parallelism of the $\psihat$ quantiles suggests that a plausible  candidate for a non-informative prior could be an improper distribution that is uniform over all $\psi$. Back in the original $\xi$ parameterisation, this corresponds to the improper prior distribution $\Pi_0(\xi)  = 1/\sqrt{1 + \xi^2}$.

This proposal for a reference prior has some attractive features. For scale-invariant problems, the $1/\sigma$ prior clearly falls away improperly as $1/\sigma$
at large $\sigma$, and for large $| \xi |$ the function $1/\sqrt{1 + \xi^2}$ falls away similarly.
Tail integrals of improper distributions naturally diverge, but when fed through the Bayesian machinery, the resulting posterior distributions can have proper, well-behaved tails (see \cite{McRobieEEE}). Earlier investigations had indicated that an improper prior uniform on $\xi$ leads to an improper posterior for the $N=3$ GPD case. However an improper prior uniform on $\psi$ (with its $~ 1/|\xi|$ behaviour at $|\xi|$ large) solves this problem, and does so without creating any awkward singularities at $\xi = 0$. Posteriors are proper, even for $N=3$.

Vertical (constant $\psi$) and horizontal (constant $\psihat$) slices through the likelihood function are shown in Figure~\ref{psigivenpsihat}. These were generated by selecting $8 \times 10^6$ values of $\psi$ distributed uniformly on $[-4,4]$.
At each value, an $N=20$ GPD sample was generated and a curve-fit estimate $\psihat$ obtained. The resulting set $\lbrace \psi, \psihat \rbrace$ was then partitioned into vertical and horizontal slices, each of width $\Delta \psi = \Delta \psihat = 0.1$, and a probability density function was fitted to the data in each slice using the Matlab {\it ksdensity} function. These sections through the likelihood function were not employed in the analysis that follows, and are presented here merely to illustrate that the $\psi = \asinh(\xi)$ re-parameterisation is a somewhat natural one that leads to well-behaved distributions.

Despite the fact that the quantile contours are strikingly almost-parallel across the full range of $\psi$, there is some small asymmetry. Closer inspection of the densities $p(\psihat|\psi)$ (Figure~\ref{psigivenpsihat}) shows that for $|\psi|$ large, the density of $\psihat$ is skewed a little to the left or right, resembling suitably-scaled versions of the scarp-and-dip shaped functions $e^{\pm x} \exp(-e^{\pm x})$. In the central region of moderate $|\xi|$, the distribution transitions between these two skewed extremes. The initial hope - that by suitable transformation, the tail-parameter analysis could be transformed into a location-parameter problem - was thus not borne out in full. However, the transformed problem is nevertheless approximately location-like. The gist of the argument is thus that, if a Bayesian is willing to countenance use of an improper uniform as the appropriate noninformative prior on a location-only problem, then it might not be too much of a stretch to countenance use of such a prior for the tail problem once suitably transformed into its almost location-like manifestation.

In summary, then, the improper distribution uniform over all $\psi$ appears to be a plausible contender for a noninformative prior, and is thus adopted as such for the rest of the paper.

\section{Constructing the Bayes-like predictor}
Having reduced the problem space from the full $\lbrace \bX, \bTheta, x_{next} \rbrace$ to the 3-dimensional $\lbrace \psi, \psihat, u_{next} \rbrace$ approximation, we now construct the Bayes-like predictor in this reduced setting. This requires constructing the predictive density (given in Equation (2) earlier) over the $\lbrace \psi, \psihat, u_{next} \rbrace$ space
\begin{equation}
p(\psihat, \psi, u_{next})  = p(u_{next} | \psihat, \psi) p(\psihat, \psi) \nonumber
\end{equation}

Although not in theory impossible, it would be extremely complicated to determine this function analytically. It is, however, rather simple to simulate the function numerically. Using the uniform prior on $\psi$, a large number of values $\psi_i$ are sampled over a wide range ($\psi_i \in [-4,4]$). At each $\psi_i$, a $N=20$ GPD is generated and a curve-fit estimate $\psihat_i$ is determined. A further singleton $x_{next}$ is then sampled from the GPD at $\psi_i$, and this is normalised using the $x_{10}$ and $x_{20}$ values of the original sample to obtain a normalised ``next sample'' $u_{next,i}$. We thus obtain a set of points $\lbrace \psi_i, \psihat_i, u_{next,i} \rbrace$ in our reduced problem space. Taking Bayes-like horizontal slices of width $\Delta \psihat = 0.1$ around a fixed estimate $\psihat$, we obtain an approximation for the density of next values $u_{next}$ given the estimate $\psihat$. A tail cdf can then be fitted to these points (using the Matlab {\it ksdensity} function) and the $T$-level quantiles $u_{T,pred}$ can be extracted. These are the Bayes-like level-$T$ predictions.

\begin{figure} [h!] \centering
   \includegraphics[width=90mm, keepaspectratio]{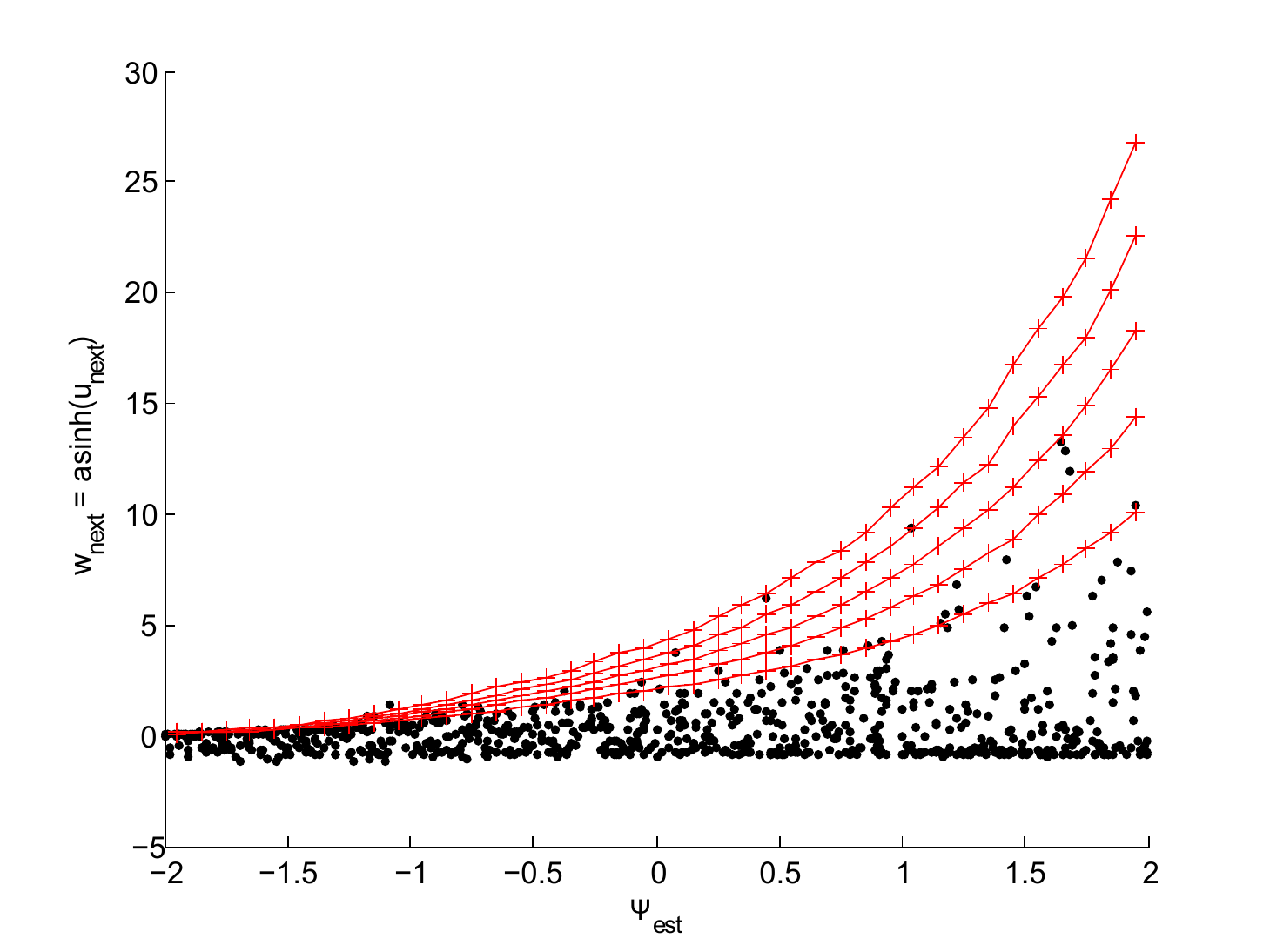}
 \caption{A small fraction of the 3D point cloud (.) projected onto the $\psihat$-$w_{next}$ plane (where $w_{next} = \asinh(u_{next})$. The level-$T$ predictions computed at each slice of $\psihat$ are superimposed in red (+) for $T = [21,50,100,200,400]$.
 }
\label{unextversusPsihat}
\end{figure}

The procedure is illustrated in Figure~\ref{unextversusPsihat}, where the 3D cloud of $\lbrace \psi_i, \psihat_i, u_{next,i} \rbrace$ points has been projected down onto the $\psihat$-$w_{next}$ plane (where $w_{next}= \asinh(u_{next})$). The cloud consisted of $8\times 10^6$ points, but to reduce visual clutter, only a tiny fraction are plotted here. The level-$T$ predictions computed at each slice (for $T = [21,50,100,200,400]$) are superimposed (red, +).

\section{Comparison of Bayes-like and Frequentist predictions}

\begin{figure} [ht!] \centering
   \includegraphics[width=70mm, keepaspectratio]{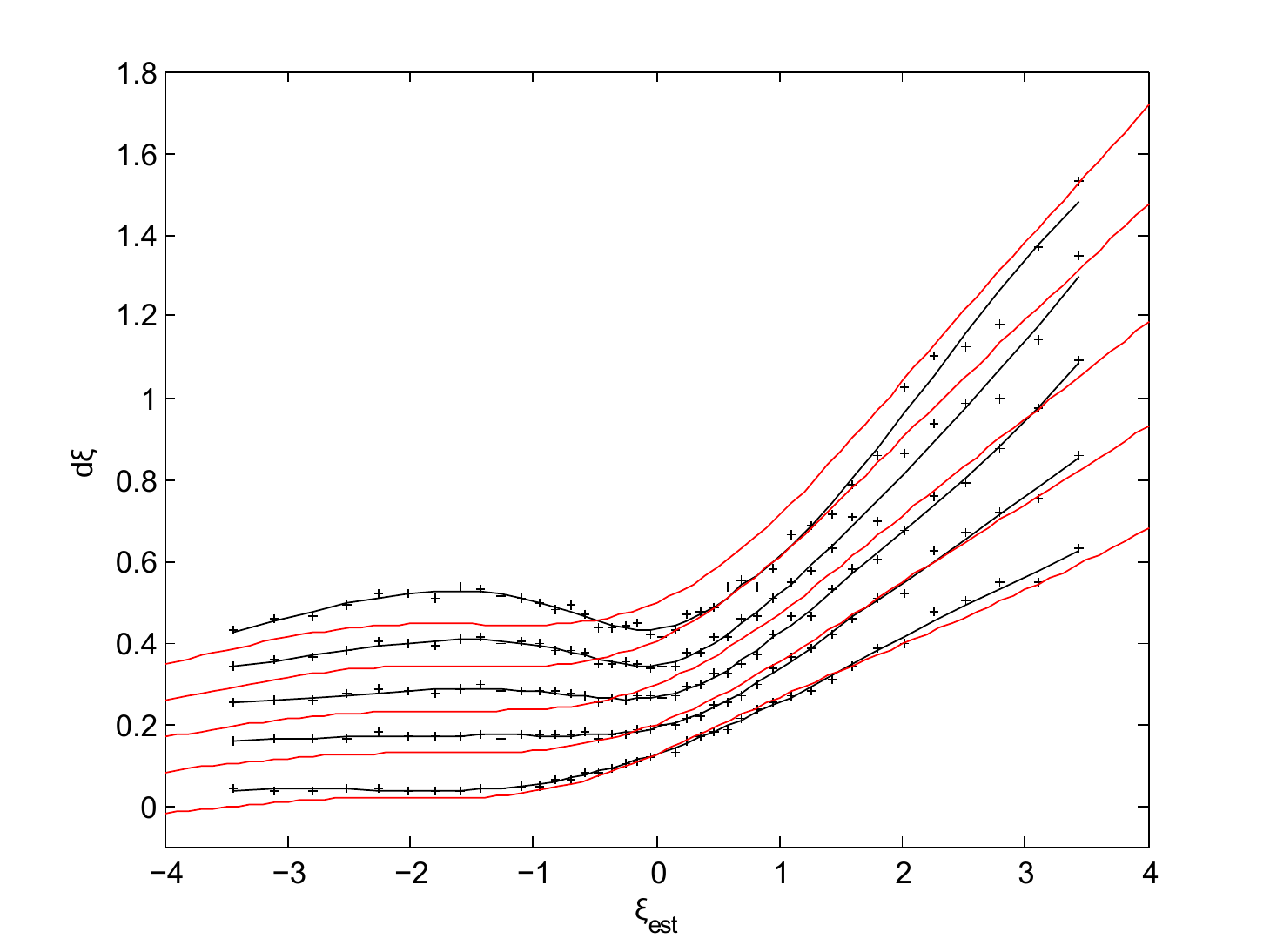}
   \includegraphics[width=70mm, keepaspectratio]{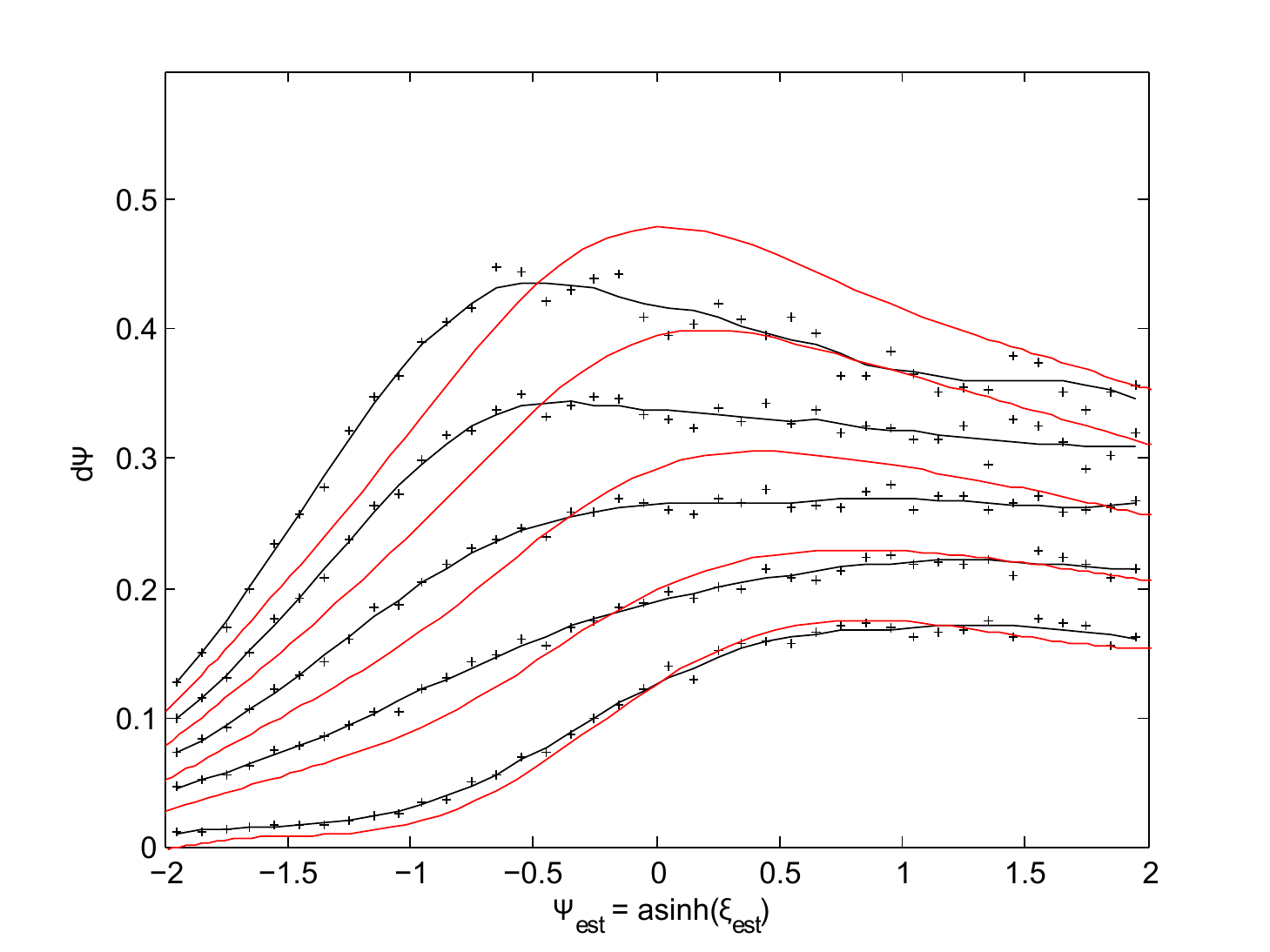}
 \caption{The computed tail parameter increments for the Bayes-like predictor (black) and the original Frequentist predictor (red).
The left hand figure shows the predictor expressed as an increment of $\xihat$ and the right-hand figure as an increment of $\psihat$.}
\label{increments}
\end{figure}

\cite{McRobiecurvefit} described the Frequentist predictor that was designed to give good probability preservation for samples of size $N=20$ drawn from a GPD of any fixed but unknown tail parameter $\xi$. The previous section has constructed an analogous Bayes-like predictor for the same problem.
The Frequentist predictor was expressed in terms of the analytical approximation underlying the curve-fit. That is, having made an estimate $\xihat$ by fitting a curve of the form
\begin{equation}
u_{i} = \frac{(G_{10}/G_i)^{\xihat}-1}{ 1 - (G_{10}/G_{20})^{\xihat}}
\end{equation}
to the data (with $G_j = (j-0.5)/N$), a level-$T$ prediction $u_{T,pred}$ could be expressed as
\begin{equation}
u_{T,pred} = \frac{(G_{10}T)^{\xi_p}-1}{ 1 - (G_{10}/G_{20})^{\xi_p}}
\end{equation}
(with $G_j = j/(N+1)$) for some $\xi_p$. That is, the prediction can be expressed as the analytical $T$-level approximation - not at the estimated value $\hat{\xi}$ - but at a higher value of the tail parameter
\begin{equation}
\xi_p = \xihat + d\xi_T
\end{equation}
Via a series of rather ad hoc function adjustments, \cite{McRobiecurvefit} arrived at functions for the increments $d\xi_T$ that led to predictors with good probability preservation. The Bayes-like predictor that has now been constructed can also be presented in this fashion, as a $T$-dependent increment to the analytical approximation at the curve-fit estimate.

The Bayes-like and Frequentist increments are compared in Figure~\ref{increments}, plotted as increments of $\xihat$ (left), and equivalently as increments of $\psihat$ (right). It can be seen that the Bayes-like and Frequentist approaches do not lead to identical results. However, the form of the solutions are decidedly similar, in terms of both the magnitude of the increment and the overall form of the increment functions. This latter point is particularly apposite. To create the Frequentist increment, the author needed to manipulate the increment function in a comparatively ad hoc manner, adding bump functions and the like, until a functional form was arrived at which delivered good probability preservation over the whole range of $\xi$. The Bayes-like increments, in contrast, were arrived at purely algorithmically, with no need for adjustments by the author. The fact that both methods arrive at the same general undulating shape of increment function adds credence to both.

\begin{figure} [ht!] \centering
   \includegraphics[width=150mm, keepaspectratio]{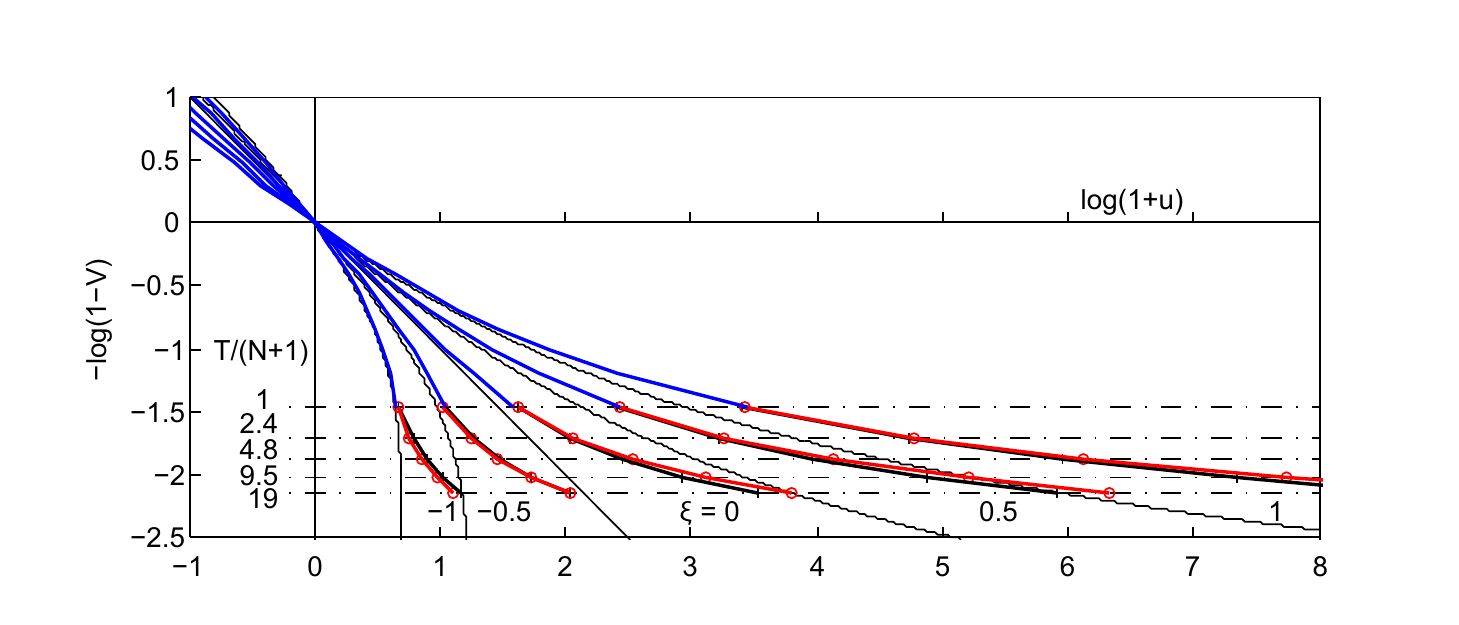}
 \caption{The extrapolation curves using $\xi_{p} = \hat{\xi} + d\xi$ (for $ \hat{\xi} = -1, -0.5, 0, 0.5 ,1$) are plotted in thick black (Bayes-like) and red (Frequentist) at the lower part of the diagram. The diagram assumes the model is GPD with $N=20$, and extrapolations are shown to recurrence levels $T$ ranging from 21 to 400, corresponding to extrapolation ratios $E_R$ from 1 to 19. In the upper part of the diagram, the data averages are shown in blue, using $G_i = i/(N+1)$ plotting positions. In the background, the basic curve-fit approximations are shown as thin black. The $y$-axis is a measure of return level, with $V = -\log(10T/21)/\log(2)$. }
\label{finalprediction}
\end{figure}

How small the differences are between the two approaches can be seen more clearly in Figure~\ref{finalprediction} where the predictions are shown as extrapolations on the basic curve-fit construction for the cases $\xi = [-1,-0.5,0,0.5,1]$. The Frequentist (red) and Bayes-like (black) extrapolations are very close, and almost overlie each other.

\begin{figure} [ht!] \centering
   \includegraphics[width=70mm, keepaspectratio]{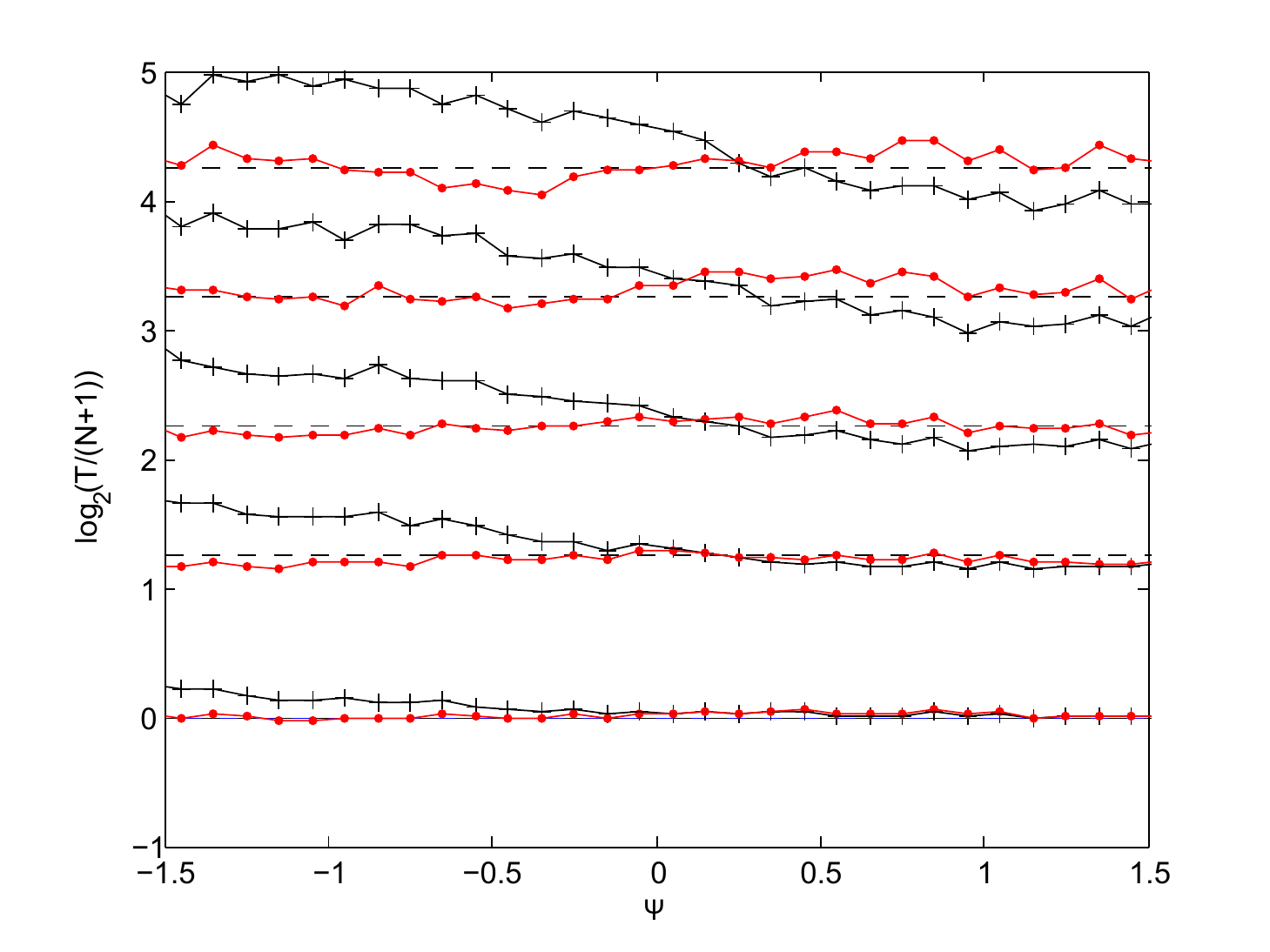}
   \includegraphics[width=70mm, keepaspectratio]{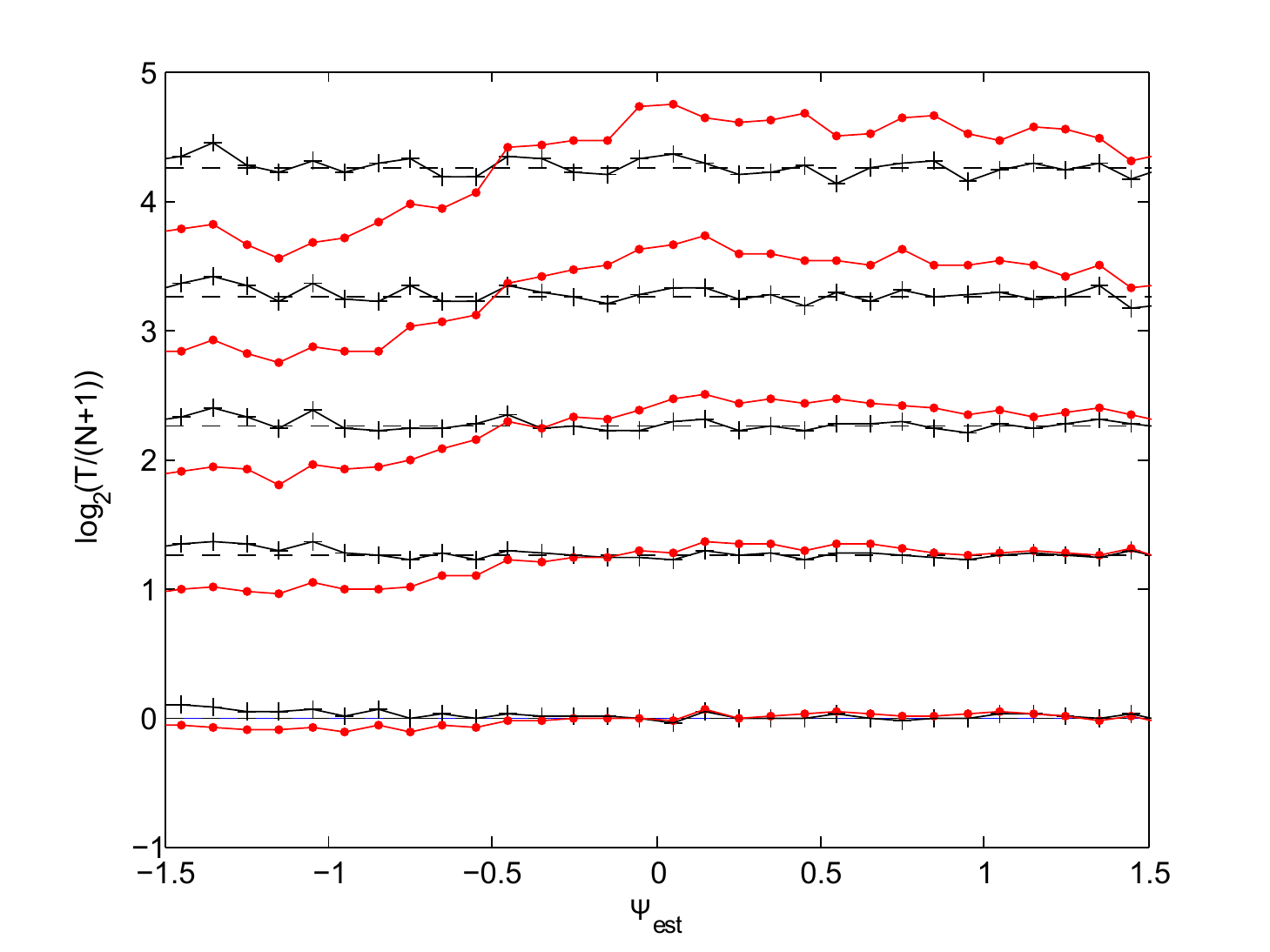}
 \caption{The delivered probability performance of the  Bayes-like predictor (black) and the original Frequentist predictor (red).
The left hand figure shows how performance varies at different tail parameters $\psi$. The Frequentist predictor gives good performance across all $\psi$, whilst the Bayes-like predictor, although not substantially different, over/underpredicts a little at $\psi$ negative/positive.
The right-hand figure is a more unusual way of considering probability preservation. Here, the probability performance is plotted with respect to the estimate $\psihat$. In this diagram, the roles are reversed. The Bayes-like predictor gives a good match to the desired return level (as it should by construction), whilst the Frequentist predictor, although close, is less precise in its delivery of the desired level. }
\label{probperform}
\end{figure}

There is no fundamental requirement for the two disparate philosophical approaches to lead to identical results. Indeed, the reason why the approaches may lead to different answers is illustrated in Figure~\ref{probperform}. There, each predictor has been applied to each point in an $8\times 10^6$ data cloud $\lbrace \psi_i, \psihat_i, u_{next,i} \rbrace$, and the number of times that the next singleton drawn $u_{next}$ exceeded each level-$T$ prediction was recorded. The left-hand diagram shows the resulting probability performance plotted against the unknown tail parameter $\xi$ (plotted in its re-parameterised form $\psi$). The Frequentist predictor, by design, delivers good probability performance across the
full range of $\psi$ considered (corresponding to $-2.1< \xi < 2.1$), whilst the Bayes-like predictor, although not too distant from delivering good probability performance, has a small tendency to over/under-predict at $\xi$ negative/positive.

The right-hand diagram in Figure~\ref{probperform} is a more unusual way to look at probability preservation. Here the graphs show what proportion of predictions were exceeded by subsequent samples of the full point cloud that lay in {\it horizontal} slices of thickness $\Delta \psihat = 0.1$ at fixed $\psihat$. Here the roles are reversed. Both approaches give good performance, but the Bayes-like predictor is almost perfect in its delivery of predictions having the desired exceedance level, whilst the Frequentist performance strays a little.  This is, of course, no surprise, because the Bayes-like predictor was constructed on horizontal slices to have exactly this property. This latter perspective, however, of looking at probability preservation with respect to estimated rather than actual parameters is really only of academic interest. Although it makes evident why the two constructions lead to slightly different results (in that taking integrals over Bayes-like horizontal slices does not lead to exactly the same predictions as taking integrals over Frequentist vertical slices), it is the performance on Frequentist vertical slices that is of practical interest.

\section{Conclusions}
The objective of this paper (and the previous paper \cite{McRobiecurvefit}) was to see how far the correspondence between the Frequentist and the $1/\sigma$ Bayesian approaches to prediction for $(\mu, \sigma)$ location and scale problems could be generalised to the three parameter $(\mu, \sigma, \xi)$ case of the GPD. The conclusion is that Frequentist and Bayes-like approaches can be constructed which have close, but not perfect, correspondence.

Given the numerous approximations involved, there may exist constructions with even closer correspondence. However, given that Figure~\ref{finalprediction} shows the predictions of the two approaches to be substantially similar, we conclude that the Frequentist predictor described in \cite{McRobiecurvefit} has a plausible claim to being a rational and prudent method of extrapolating to extremes beyond the span of historical data in the case of vague prior information of parameters.

The Bayes-like predictor described in this paper was constructed algorithmically, and can be readily generalized to sample sizes $N$ other than 20. Although the use of a reference prior removes subjectivity in the pure GPD case, when predictions are to be made using tails of non-GPD data, judgement will still be required in choosing the number $N$ of upper order statistics to include, for both Frequentist and Bayes-like approaches.

\bibliography{bayeslikebib}

\end{document}